\documentclass[12pt]{amsart}

\usepackage{amsmath,amssymb,amsthm,amscd,amsfonts}
\usepackage{verbatim}
\usepackage{comment}
\usepackage{multirow}
\usepackage{mathtools}
\usepackage{enumitem}
\usepackage{pgf,tikz}
\usepackage{tikz-cd}
\usetikzlibrary{arrows}
\usetikzlibrary{cd} 
\usepackage[normalem]{ulem}
\usepackage{graphicx}
\usepackage{wasysym}
\usepackage{csquotes}
\usepackage{bbm}
\usepackage{bm}
\usepackage{color}
\usepackage{tabularx}
\usepackage{stmaryrd}
\usepackage[margin=1in]{geometry} 

\numberwithin{equation}{section}

\usepackage{pdfpages}
\usepackage{hyperref}

\newcommand{\Ext}{\operatorname{Ext}}
\newcommand{\GL}{\operatorname{GL}}
\newcommand{\Sym}{\operatorname{Sym}}

\newcommand{\ZZ}{\mathbb{Z}}
\def\PP{{\mathbf P}}

\newcommand{\OO}{\mathcal{O}}
\newcommand{\RR}{\mathcal{R}}
\newcommand{\oo}{\otimes}
\newcommand{\bb}[1]{\mathbb{#1}}

\newcommand{\mc}[1]{\mathcal{#1}}

\def\lra{\longrightarrow}

\def\kk{{\mathbf k}}

\newcommand{\op}[1]{\operatorname{#1}}
\newcommand{\extt}[1]{\small\texttt{#1}}

\DeclareMathOperator{\HH}{H}

\DeclareMathOperator{\ch}{char}

\theoremstyle{plain}

\newtheorem{thm}{Theorem}[section]

\newtheorem{con}[thm]{Conjecture}

\theoremstyle{definition}

\newtheorem{ex}[thm]{Example}
\theoremstyle{remark}

\begin{document}

\title[Computing cohomology on the incidence correspondence]{Computing the cohomology of line bundles on the incidence correspondence and related invariants}

\author{Annet Kyomuhangi}
\address{Department of Mathematics, Busitema University, P.O. Box 236, Tororo}
\email{annet.kyomuhangi@gmail.com}

\author{Emanuela Marangone}
\address{Department of Mathematics, University of Manitoba, Winnipeg, MB R3T2M8 \newline
\indent PIMS - Pacific Institute for the Mathematical Sciences}
\email{emanuela.marangone@umanitoba.ca}

\author{Claudiu Raicu}
\address{Department of Mathematics, University of Notre Dame, 255 Hurley, Notre Dame, IN 46556\newline
\indent Institute of Mathematics ``Simion Stoilow'' of the Romanian Academy}
\email{craicu@nd.edu}

\author{Ethan Reed}
\address{Department of Mathematics, University of Notre Dame, 255 Hurley, Notre Dame, IN 46556}
\email{ereed4@nd.edu}


\date{\today}

\begin{abstract}
We describe the package \texttt{IncidenceCorrespondenceCohomology} for the computer algebra system \textit{Macaulay2}. The main feature concerns the computation of characters and dimensions for the cohomology groups of line bundles on the incidence correspondence (the partial flag variety parametrizing pairs consisting of a point in projective space and a hyperplane containing it). Additionally, the package provides tools for  (1) computing the multiplication in the graded Han--Monsky representation ring, (2) determining the splitting type of vector bundles of principal parts on the projective line, and (3) testing the weak and strong Lefschetz properties for Artinian monomial complete intersections.
\end{abstract}

\maketitle

\section{Introduction}
The Borel--Weil--Bott theorem \cite{borel-weil,bott} is a fundamental result in algebraic geometry and representation theory, describing the cohomology groups of line bundles on flag varieties over a field of characteristic zero. The corresponding question over fields of positive characteristic remains open, and has seen renewed interest in recent years (see \cite{rai-vdb,GRV} for some recent results and conjectures, and for references to some of the classical work on this topic). An interesting special case where the cohomology groups are highly sensitive to the characteristic is that of the \textbf{incidence correspondence} -- the partial flag variety parametrizing pairs consisting of a point in projective space and a hyperplane containing it. In \cite{KMRR} we obtain general character and dimension formulas (mostly recursive) for the cohomology of line bundles on the incidence correspondence in arbitrary characteristic, and we point out surprising connections with several other topics:
\begin{itemize}
    \item the multiplication in the graded version of the Han--Monsky representation ring,
    \item the splitting type of vector bundles of principal parts on the projective line,
    \item the weak Lefschetz property for Artinian monomial complete intersections.
\end{itemize}
The current article can be read as a companion to \cite{KMRR}. It describes the package \texttt{IncidenceCorrespondenceCohomology} for the computer algebra system \textit{Macaulay2}, where we implement algorithms to compute cohomology and related invariants based on the theoretical results from \cite{KMRR}.

\smallskip

\noindent{\bf Organization.} In Section~\ref{sec:cohIC} we discuss the calculation of cohomology of line bundles on the incidence correspondence, and the closely related problem of determining the cohomology for twists of divided powers of the cotangent sheaf on projective space. In Section~\ref{sec:split-pparts} we describe the calculation of the splitting type of vector bundles of principal parts on the projective line. In Section~\ref{sec:han-monsky} we discuss the multiplication in the graded Han--Monsky representation ring, with applications to testing the weak and strong Lefschetz properties for Artinian monomial complete intersections. In addition, we implement methods to test the weak Lefschetz property for more general Artinian algebras.

\section{Cohomology of line bundles on the incidence correspondence}\label{sec:cohIC}

Let $n\geq 2$ and consider the projective space $\PP=\mathbb{P}(\kk^n)=\bb{P}^{n-1}$, where $\kk$ is an algebraically closed field of characteristic $p$. We write $X$ for the \textbf{incidence correspondence} defined by
$$X=\{(x,H)| \ x\in H \}\subset \PP\times \PP^\vee,$$
where $\PP^\vee$ is the dual projective space parametrizing hyperplanes in $\PP$. Each line bundle on~$X$ is the restriction of a line bundle on $\PP\times\PP^{\vee}$, so it has the form 
$$\OO_X(a,b) = \OO_{\PP \times \PP^{\vee}}(a,b)_{|_X}\quad\text{ for }(a,b)\in\ZZ^2.$$
The cohomology groups of $\OO_X(a,b)$ are finite dimensional representations of $\GL_n$. For any such representation $W$, we consider its character
\[[W] = \sum_{(i_1,\cdots,i_n)\in\bb{Z}^n} \dim(W_{(i_1,\cdots,i_n)})\cdot z_1^{i_1}\cdots z_n^{i_n} \]
which is the symmetric Laurent monomial in $\bb{Z}[z_1^{\pm 1},\cdots,z_n^{\pm 1}]$ that records the eigenspace decomposition of $W$ relative to the action of the diagonal torus $(\kk^\times)^n$ in $\GL_n$. The dimension of $W$ can be recovered from its character by setting $z_i=1$ for each $i$. The main goal of our package is to provide a computational tool for the characters
\begin{equation}\label{eq:def-hj-OXab} 
h^i(\OO_X(a,b)) = \left[\HH^i\left(X,\mc{O}_X(a,b)\right)\right],\quad\text{ where }i\geq 0,\text{ and }a,b\in\bb{Z}.
\end{equation}
In characteristic zero this is a simple application of the Borel--Weil--Bott theorem, and for certain values of the parameters (such as $a,b\geq 0$) the characters \eqref{eq:def-hj-OXab} are independent of the characteristic (see the discussion following \cite[Theorem~1.1]{gao-raicu} and Section~\ref{subsec:inc-coh} below).

A key insight in \cite{gao-raicu,KMRR} is that the more subtle character calculations which depend on the characteristic (such as when $a\geq 0$, $b\leq -n$) can be performed more effectively as cohomology calculations for some higher rank vector bundles on $\PP$. More precisely, we let $\Omega$ denote the cotangent sheaf on $\PP$, let $\mc{R}=\Omega(1)$ denote the universal rank $(n-1)$ subsheaf, let $D^d\mc{R}$ denote the $d$-th divided power of $\mc{R}$, and consider the characters
\begin{equation}\label{eq:def-hj-Ddre} 
 h^i(D^d\mc{R}(e)) = \left[\HH^i(\PP,D^d\mc{R}(e)) \right],\quad\text{ where }d,i\geq 0,\text{ and }e\in\bb{Z}.
\end{equation}
Using \cite[(2.6)]{KMRR} (see also \cite[(2.12)]{gao-raicu}), we have
\begin{equation}\label{div-inc}
    h^i(D^d\mc{R}(e)) \cdot z_1\cdots z_n =  h^{i+n-2}\left(\OO_X(e+1,-d-n+1)\right)\quad\text{ for }d,i\geq 0,\text{ and }e\in\bb{Z}.
\end{equation}
Following \cite{KMRR}, in Section~\ref{subsec:rec-coh-divpows} we discuss the recursive calculation of the characters~\eqref{eq:def-hj-Ddre}, and in Section~\ref{subsec:nim-coh-divpows} we describe a non-recursive calculation in characteristic $p=2$.

\subsection{Cohomology for divided powers of the universal subsheaf}\label{subsec:rec-coh-divpows}

In this section we discuss a recursive method for computing the characters \eqref{eq:def-hj-Ddre} and the corresponding dimensions when $d\geq 0$ and $e\geq -1$. Under this restriction on the parameters, one has that $h^i(D^d\mc{R}(e))=0$ for $i\neq 0,1$. Moreover, one has
\begin{equation}\label{eq:h1-vs-h0}
h^i(D^d\mc{R}(e)) = h^{1-i}(D^{e+1}\mc{R}(d-1))\quad\text{ for }i=0,1,
\end{equation}
which reduces the calculation to the situation when $e\geq d-1$. The base step of the recursion is when $d<p$ where the cohomology is computed using \cite[(1.2)]{KMRR}, while the recursive step is given by \cite[Theorem~1.1]{KMRR}.

The method we implement for this calculation is called \texttt{recursiveDividedCohomology}, and it has an optional input \texttt{FindCharacter}. The default value of \texttt{FindCharacter} is \texttt{false} in which case our function returns the dimension of cohomology, while in the case when \texttt{FindCharacter} is set to \texttt{true}, the function outputs the cohomology character.

The main instance of this function takes input a list $\{ i, p, d, e, n\}$, where $i\in\{0,1\}$, $d\geq 0$, $e\geq -1$, and it outputs the character or dimension of the cohomology group $\HH^i(\bb{P}(\kk^n),D^d\mc{R}(e))$ when the field $\kk$ has characteristic $p$.

\begin{ex} The following computes $\HH^1\left(\mathbb{P}^3, D^3\RR(2)\right)$ in characteristic $p=3$:

\medskip

\noindent\extt{i1: loadPackage "IncidenceCorrespondenceCohomology"\\
o1 = IncidenceCorrespondenceCohomology\\
o1: Package\\
i2: recursiveDividedCohomology(\{1,3,3,2,4\})\\
o2 = 16\\
i3: recursiveDividedCohomology(\{1,3,3,2,4\}, FindCharacter => true)\\
o3 = $\mathtt{x_{1}^{2}x_{2}^{2}x_{3}+x_{1}^{2}x_{2}x_{3}^{2}+x_{1}x_{2}^{2}x_{3}^{2}+x_{1}^{2}x_{2}^{2}x_{4}+x_{1}^{2}x_{2}x_{3}x_{4}+x_{1}x_{2}^{2}x_{3}x_{4}+x_{1}^{2}x_{3}^{2}x_{4}+x_{1}x_{2}x_{3}^{2}x_{4}+}$\\ \phantom{o3 = }$\mathtt{x_{2
     }^{2}x_{3}^{2}x_{4}+x_{1}^{2}x_{2}x_{4}^{2}+x_{1}x_{2}^{2}x_{4}^{2}+x_{1}^{2}x_{3}x_{4}^{2}+x_{1}x_{2}x_{3}x_{4}^{2}+x_{2}^{2}x_{3}x_{4}^{2}+x_{1}x_{3}^{2}x_{4}^{2}+x_{2}x_{3}^{2}x_{4}^{2}}$\\
o3: ZZ[$\mathtt{x_1,\cdots, x_4}$]}
\end{ex}

\begin{ex}\label{ex:h1D4R2} The following computes $\HH^1\left(\mathbb{P}^2, D^4\RR(2)\right)$ in characteristic $p=3$:
\medskip

\noindent\extt{i4:  recursiveDividedCohomology(\{1,3,4,2,3\}, FindCharacter => true)\\
o4 = $\mathtt{x_{1}^{3}x_{2}^{3}+x_{1}^{3}x_{2}^{2}x_{3}+x_{1}^{2}x_{2}^{3}x_{3}+x_{1}^{3}x_{2}x_{3}^{2}+2\,x_{1}^{2}x_{2}^{2}x_{3}^{2}+x_{1}x_{2}^{3}x_{3}^{2}+x_{1}^{3}x_{3}^{3}+x_{1}^{2}x_{2}x_{3}^{3}+x_{1}x_{2
      }^{2}x_{3}^{3}+x_{2}^{3}x_{3}^{3}}$\\
o4: ZZ[$\mathtt{x_1,\dots, x_3}$]}
\smallskip

\noindent which is equivalent to computing $\HH^0\left(\mathbb{P}^2, D^3\RR(3)\right)$:
\smallskip

\noindent\extt{i5:  recursiveDividedCohomology(\{0,3,3,3,3\}, FindCharacter => true)\\   
o5 = $\mathtt{x_{1}^{3}x_{2}^{3}+x_{1}^{3}x_{2}^{2}x_{3}+x_{1}^{2}x_{2}^{3}x_{3}+x_{1}^{3}x_{2}x_{3}^{2}+2\,x_{1}^{2}x_{2}^{2}x_{3}^{2}+x_{1}x_{2}^{3}x_{3}^{2}+x_{1}^{3}x_{3}^{3}+x_{1}^{2}x_{2}x_{3}^{3}+x_{1}x_{2
      }^{2}x_{3}^{3}+x_{2}^{3}x_{3}^{3}}$\\
o5: ZZ[$\mathtt{x_1,\dots, x_3}$]\\
i6: recursiveDividedCohomology(\{1,3,4,2,3\})\\
o6 = 11\\
i7: recursiveDividedCohomology(\{1,3,4,2,3\})==recursiveDividedCohomology(\{0,3,3,3,3\})\\
o7 = true}
\end{ex}

An alternative instance of the method \texttt{recursiveDividedCohomology} allows one to specify the character ring $R$ instead of the number of variables: the inputs are then a list $\{ i, p, d, e\}$ where the entries have the same significance as before, and a ring $R$. The output is then the cohomology character expressed as an element of $R$, and the optional input \texttt{FindCharacter} is no longer relevant. Under our working assumptions $d\geq 0$, $e\geq -1$, the cohomology characters \eqref{eq:def-hj-Ddre} are usual polynomials (as opposed to Laurent polynomials) so it is not necessary for the variables in $R$ to have inverses.

\begin{ex}
As noted in Example~\eqref{ex:h1D4R2} one has $h^1(D^4\RR(2))=h^0(D^3\RR(3))$, and this can be validated in Macaulay2 as long as the two characters belong to the same ring:
\smallskip

\noindent\extt{i8: R=ZZ[$\mathtt{z_1, z_2, z_3}$]\\
o8 = R\\
o8: PolynomialRing\\
i9: recursiveDividedCohomology(\{1,3,4,2\},R)\\
o9 = $\mathtt{z_{1}^{3}z_{2}^{3}+z_{1}^{3}z_{2}^{2}z_{3}+z_{1}^{2}z_{2}^{3}z_{3}+z_{1}^{3}z_{2}z_{3}^{2}+2\,z_{1}^{2}z_{2}^{2}z_{3}^{2}+z_{1}z_{2}^{3}z_{3}^{2}+z_{1}^{3}z_{3}^{3}+z_{1}^{2}z_{2}z_{3}^{3}+z_{1}z_{2
      }^{2}z_{3}^{3}+z_{2}^{3}z_{3}^{3}}$\\
o9: R\\
i10: recursiveDividedCohomology(\{1,3,4,2\},R)==recursiveDividedCohomology(\{0,3,3,3\},R)
o10 = true}    
\end{ex}

\subsection{Non-recursive cohomology in characteristic \boldmath$2$}\label{subsec:nim-coh-divpows}

We now restrict to characteristic $p=2$ and discuss a non-recursive method to compute the characters \eqref{eq:def-hj-Ddre}, still under the assumption $d\geq 0$, $e\geq -1$. As explained earlier, \eqref{eq:h1-vs-h0} allows one to further assume $e\geq d-1$, in which case we can apply \cite[Theorem~1.2]{KMRR} in order to compute $h^1(D^d\mc{R}(e))$. Since the cohomology is concentrated in degrees $0,1$, in order to compute $h^0(D^d\mc{R}(e))$ it then suffices to know the Euler characteristic of $D^d\mc{R}(e)$. This is independent of the characteristic of $\kk$ and is given by
\[h^0(D^d\mc{R}(e)) - h^1(D^d\mc{R}(e)) = h_d\cdot h_e - h_{d-1}\cdot h_{e+1},\]
where $h_d$ denotes the degree $d$ \textbf{complete symmetric polynomial}, and is computed by the formula
\begin{equation}\label{eq:def-hd}
    h_d = \sum_{i_1+\cdots+i_n=d} z_1^{i_1}\cdots z_n^{i_n}
\end{equation}

The method we implement for this calculation is called \texttt{nimDividedCohomology}, and takes an optional input \texttt{FindCharacter} with the same significance and default value as in Section~\ref{subsec:rec-coh-divpows}. The main instance of the method \texttt{nimDividedCohomology} takes inputs $i,d,e,n$, where $i\in\{0,1\}$, $d\geq 0$, $e\geq -1$, and it outputs the character or dimension of the cohomology group $\HH^i(\bb{P}(\kk^n),D^d\mc{R}(e))$ when the field $\kk$ has characteristic $p=2$.

\begin{ex}\label{ex:H1D3R3-char2} The following computes $\HH^1\left(\mathbb{P}^4, D^3\RR(3)\right)$ in characteristic $2$:
\smallskip

\noindent\extt{i11: nimDividedCohomology(1,3,3,5, FindCharacter => true)\\
o11 = $\mathtt{x_{1}^{2}x_{2}x_{3}x_{4}x_{5}+x_{1}x_{2}^{2}x_{3}x_{4}x_{5}+x_{1}x_{2}x_{3}^{2}x_{4}x_{5}+x_{1}x_{2}x_{3}x_{4}^{2}x_{5}+x_{1}x_{2}x_{3}x_{4}x_{5}^{2}}$\\
o11: ZZ[$\mathtt{x_1,\dots, x_3}$]}
\smallskip

\noindent which can also be computed using the method from Section~\ref{subsec:rec-coh-divpows}:
\smallskip

\noindent\extt{i12: recursiveDividedCohomology(\{1,2,3,3,5\}, FindCharacter => true)  \\
o12 = $\mathtt{x_{1}^{2}x_{2}x_{3}x_{4}x_{5}+x_{1}x_{2}^{2}x_{3}x_{4}x_{5}+x_{1}x_{2}x_{3}^{2}x_{4}x_{5}+x_{1}x_{2}x_{3}x_{4}^{2}x_{5}+x_{1}x_{2}x_{3}x_{4}x_{5}^{2}}$\\
o12: ZZ[$\mathtt{x_1,\dots, x_3}$]\\
i13: recursiveDividedCohomology(\{1,2,3,3,5\})\\
o13 = 5\\
i14: recursiveDividedCohomology(\{1,2,3,3,5\})==nimDividedCohomology(1,3,3,5)\\
o14 = true}
\end{ex}

As it was the case for the method \texttt{recursiveDividedCohomology} in Section~\ref{subsec:rec-coh-divpows}, an alternative instance of \texttt{nimDividedCohomology} allows one to specify a character ring $R$ instead of the number of variables $n$. In this case the output is necessarily a character polynomial which belongs to $R$, and the option \texttt{FindCharacter} is no longer relevant.

\begin{ex} The following computes $\HH^1\left(\mathbb{P}^3, D^4\RR(7)\right)$ in characteristic $2$:
\smallskip

\noindent\extt{i15: R=ZZ[$\mathtt{z_1..z_4}$]\\
o15 = R\\
o15: PolynomialRing\\
i16: nimDividedCohomology(1,4,7,R)\\
o16: $\mathtt{z_{1}^{3}z_{2}^{3}z_{3}^{3}z_{4}^{2}+z_{1}^{3}z_{2}^{3}z_{3}^{2}z_{4}^{3}+z_{1}^{3}z_{2}^{2}z_{3}^{3}z_{4}^{3}+z_{1}^{2}z_{2}^{3}z_{3}^{3}z_{4}^{3}}$}
\end{ex}

While both the methods \texttt{recursiveDividedCohomology} and \texttt{nimDividedCohomology} can be used to perform the computation in characteristic $2$, the second one is often faster.
In the next examples, we show the running times of the two different methods to compute the characters, so in both cases we set the option \texttt{FindCharacter} to \texttt{true}.
\begin{ex} Time comparison for the computation of $h^1\left(D^3\RR(3)\right)$ when $n=5$ (as in Example~\ref{ex:H1D3R3-char2})

      \begin{tabularx}{0.95\textwidth} { 
  | >{\centering\arraybackslash}X 
  | >{\centering\arraybackslash}X | }
 \hline
    \texttt{recursiveDividedCohomology} &  \texttt{nimDividedCohomology} \\
 \hline
  0.0952897 seconds &     0.0038127 seconds\\
\hline
\end{tabularx}
\end{ex}
\begin{ex} Time comparison for the computation of $h^1\left(D^6\RR(9)\right)$  with $n=7$ 
\smallskip

      \begin{tabularx}{0.95\textwidth} { 
  | >{\centering\arraybackslash}X 
  | >{\centering\arraybackslash}X | }
 \hline
    \texttt{recursiveDividedCohomology} &  \texttt{nimDividedCohomology} \\
 \hline
    30.8806 seconds &    4.16985 seconds\\
\hline
\end{tabularx}
\smallskip

\noindent In this case the dimension of $\HH^1\left(\mathbb{P}^6,D^6\RR(9)\right)$ is $35637$. 
We can also notice a significant difference in times needed for the two methods to compute the dimensions: 
\smallskip

   \begin{tabularx}{0.95\textwidth} { 
  | >{\centering\arraybackslash}X 
  | >{\centering\arraybackslash}X | }
 \hline
    \texttt{recursiveDividedCohomology} &  \texttt{nimDividedCohomology} \\
 \hline
   0.0237189 seconds &    0.0016667 seconds\\
\hline
\end{tabularx}
\end{ex}

\begin{ex}  Time comparison for the computation of $h^1\left(D^{17}\RR(18)\right)$  with $n=4$ 

      \begin{tabularx}{0.95\textwidth} { 
  | >{\centering\arraybackslash}X 
  | >{\centering\arraybackslash}X | }
 \hline
    \texttt{recursiveDividedCohomology} &  \texttt{nimDividedCohomology} \\
 \hline
   20.8714 seconds &     3.1881 seconds\\
\hline
\end{tabularx}
\smallskip

\noindent In this case the dimension of $\HH^1\left(\mathbb{P}^3,D^{17}\RR(18)\right)$ is $9040$.
\end{ex}

\subsection{Computing cohomology on the incidence correspondence}\label{subsec:inc-coh}

In this section we discuss the calculation of the cohomology of line bundles on the incidence correspondence, which takes advantage of the identification \eqref{div-inc} and the recursive algorithm in Section~\ref{subsec:rec-coh-divpows}. We begin by reviewing some general results regarding the cohomology of $\mc{O}_X(a,b)$, following \cite[Introduction]{gao-raicu}.

The involution $W\mapsto W^{\vee}$ on representations of $\GL_n$ induces an involution on the character ring sending $z_i \mapsto z_i^{-1}$ for $i=1,\cdots,n$. We have
\begin{equation}\label{eq:swap-ab} 
h^i(\mc{O}_X(a,b)) = h^i(\mc{O}_X(b,a))^{\vee}.
\end{equation}
The dimension of $X$ is $2n-3$ and the canonical line bundle is $\omega_X=\mc{O}_X(1-n,1-n)$, hence 
\begin{equation}\label{eq:Serre-on-X} 
h^i(\mc{O}_X(a,b)) = h^{2n-3-i}(\mc{O}_X(1-n-a,1-n-b))^{\vee}.
\end{equation}
If $a,b\geq 0$ then $h^i(\mc{O}_X(a,b))=0$ for $i\neq 0$ and using notation \eqref{eq:def-hd} we have
\[ h^0(\mc{O}_X(a,b)) = h_a\cdot h_b^{\vee} - h_{a-1}\cdot h_{b-1}^{\vee}.\]
It follows that if $a,b\leq 1-n$ then $h^i(\mc{O}_X(a,b))=0$ for $i\neq 2n-3$, and $h^{2n-3}(\mc{O}_X(a,b))$ can be computed via \eqref{eq:Serre-on-X}. If either $2-n\leq a\leq -1$ or $2-n\leq b\leq -1$ then $h^i(\mc{O}_X(a,b))=0$ for all $i$. It follows using \eqref{eq:swap-ab} that it remains to consider the case when $a\geq 0$ and $b\leq 1-n$, when we get using \eqref{div-inc} that
\[ h^i(\mc{O}_X(a,b)) = h^{i-n+2}(D^d\mc{R}(e))\cdot z_1\cdots z_n,\text{ where }e=a-1,\ d=1-n-b.\]
Notice that in the above we have $d\geq 0$, $e\geq -1$, hence the method in Section~\ref{subsec:rec-coh-divpows} allows one to compute $h^i(\mc{O}_X(a,b))$, which moreover can only be non-zero when $i=n-2$ or $i=n-1$.

Based on the discussion above we implement the method \texttt{incidenceCohomology} to compute the dimensions and characters of cohomology groups of $\OO_X(a,b)$, where the non-trivial part of the calculation is based on the method \texttt{recursiveDividedCohomology}. As in the previous sections, the method takes the optional input \texttt{FindCharacter} with default value \texttt{false}. For the main instance of the method \texttt{incidenceCohomology}, the input is a list of integers $\{i,p,a,b,n\}$, and the output is the dimension or character of $\HH^i\left(X, \OO_X(a,b)\right)$, where $X$ is the incidence correspondence of dimension $2n-3$ over a field $\kk$ of characteristic $p$. 

\begin{ex} The following computes the dimensions of the non-zero cohomology groups of $\OO_X(5,-7)$ for $n=4$ and $p=3$:
\smallskip

\noindent\extt{i17: incidenceCohomology(\{2,3,5,-7,4\})\\
o17 = 120\\
i18: incidenceCohomology(\{3,3,5,-7,4\})\\
o18 = 15}
\end{ex}

As in the case of the methods \texttt{recursiveDividedCohomology} and \texttt{nimDividedCohomology}, an alternative instance of the method \texttt{incidenceCohomology} allows the user to input a character ring $R$ instead of the parameter $n$. As opposed to the previous situations where the characters were given by polynomials, in this case the ring $R$ must be a Laurent polynomial ring. The output is then an element of $R$ and the option \texttt{FindCharacter} is not relevant.

\begin{ex} To compute the character of $\HH^3\left(X, \OO_X(5,-7)\right)$  for $n=4$ we can use the option \texttt{FindCharacter}:
\smallskip

\noindent\extt{i19: incidenceCohomology(\{3,3,5,-7,4\}, FindCharacter => true)\\
o19 = $\mathtt{x_{1}^{4}x_{2}^{3}x_{3}^{3}x_{4}^{2}+x_{1}^{4}x_{2}^{3}x_{3}^{2}x_{4}^{3}+x_{1}^{4}x_{2}^{2}x_{3}^{3}x_{4}^{3}+x_{1}^{3}x_{2}^{4}x_{3}^{3}x_{4}^{2}+x_{1}^{3}x_{2}^{4}x_{3}^{2}x_{4}^{3}+x_{1}^{3}x_{2}^{3}x_{3}^{4}x_{4}^{2}+3x_{1}^{3}x_{2}^{3}x_{3}^{3}x_{4}^{3}+}$\\
\phantom{o19 = }$\mathtt{x_{1}^{3}x_{2}^{3}x_{3}^{2}x_{4}^{4}+x_{1}^{3}x_{2}^{2}x_{3}^{4}x_{4}^{3}+x_{1}^{3}x_{2}^{2}x_{3}^{3}x_{4}^{4}+x_{1}^{2}x_{2}^{4}x_{3}^{3}x_{4}^{3}+x_{1}^{2}x_{2}^{3}x_{3}^{4}x_{4}^{
      3}+x_{1}^{2}x_{2}^{3}x_{3}^{3}x_{4}^{4}}$\\}

\noindent or we can specify the character ring $R$ as follows:
\smallskip

\noindent\extt{i20:  R = ZZ[$\mathtt{z_1..z_4}$, Inverses=>true, MonomialOrder=>Lex]\\
o20 = R\\
o20: PolynomialRing\\
i21: incidenceCohomology(\{3,3,5,-7\},R)\\
o21 = $\mathtt{z_{1}^{4}z_{2}^{3}z_{3}^{3}z_{4}^{2}+z_{1}^{4}z_{2}^{3}z_{3}^{2}z_{4}^{3}+z_{1}^{4}z_{2}^{2}z_{3}^{3}z_{4}^{3}+z_{1}^{3}z_{2}^{4}z_{3}^{3}z_{4}^{2}+z_{1}^{3}z_{2}^{4}z_{3}^{2}z_{4}^{3}+z_{1}^{3}z_{2}^{3}z_{3}^{4}z_{4}^{2}+3z_{1}^{3}z_{2}^{3}z_{3}^{3}z_{4}^{3}+}$\\
\phantom{o21 = }$\mathtt{z_{1}^{3}z_{2}^{3}z_{3}^{2}z_{4}^{4}+z_{1}^{3}z_{2}^{2}z_{3}^{4}z_{4}^{3}+z_{1}^{3}z_{2}^{2}z_{3}^{3}z_{4}^{4}+z_{1}^{2}z_{2}^{4}z_{3}^{3}z_{4}^{3}+z_{1}^{2}z_{2}^{3}z_{3}^{4}z_{4}^{3}+z_{1}^{2}z_{2}^{3}z_{3}^{3}z_{4}^{4}}$\\}

\noindent In the ring $R$ we can verify the equation \eqref{div-inc} for this example:
\smallskip

\noindent\extt{i22: f=recursiveDividedCohomology(\{1,3,4,4\},R)\\
o22 = $\mathtt{z_{1}^{3}z_{2}^{2}z_{3}^{2}z_{4}+z_{1}^{3}z_{2}^{2}z_{3}z_{4}^{2}+z_{1}^{3}z_{2}z_{3}^{2}z_{4}^{2}+z_{1}^{2}z_{2}^{3}z_{3}^{2}z_{4}+z_{1}^{2}z_{2}^{3}z_{3}z_{4}^{2}+z_{1}^{2}z_{2}^{2}z_{3}^{3}z_{4}+3z_{1}^{2}z_{2}^{2}z_{3}^{2}z_{4}^{2}+}$\\
\phantom{o22 = }$\mathtt{z_{1}^{2}z_{2}^{2}z_{3}z_{4}^{3}+z_{1}^{2}z_{2}z_{3}^{3}z_{4}^{2}+z_{1}^{2}z_{2}z_{3}^{2}z_{4}^{3}+z_{1}z_{2}^{3}z_{3}^{2}z_{4}^{2}+z_{1}z_{2}^{2}z_{3}^{3}z_{4}^{2}+z_{1}z_{2}^{2}z_{3}^{2}z_{4}^{3}}$\\
i23: incidenceCohomology(\{3,3,5,-7\},R)==$\mathtt{f\cdot z_1z_2z_3z_4}$\\
o23 = true}
\end{ex}

\section{Splitting type of vector bundles of principal parts}\label{sec:split-pparts}

In this section we consider as before a field $\kk$ of characteristic $p$, and we focus on the projective line $\PP=\bb{P}(\kk^2)=\bb{P}^1$. Following \cite[Section~3]{KMRR} we study the splitting as a sum of line bundles for vector bundles $\mc{F}^d_r$ defined by a short exact sequence
\begin{equation}\label{eq:ses-Fdr} 0 \lra \mc{F}^d_r \lra D^d(\kk^2) \oo \OO_{\PP} \lra D^{d-r}(\kk^2) \oo \OO_{\PP}(r) \lra 0
\end{equation}
Our interest in the bundles $\mc{F}^d_r$ comes from their relation to the graded Han--Monsky representation ring discussed in Section~\ref{sec:han-monsky}, but as we explain in Section~\ref{subsec:pparts} these bundles are also closely related to vector bundles of principal parts associated to line bundles on $\PP$.

We let $T=\bb{G}_m\times \bb{G}_m$ denote the diagonal torus in $\GL_2$ and recall that the irreducible $T$-representations $L_{u,v}$ are $1$-dimensional, indexed by pairs $(u,v)\in\bb{Z}^2$. The $T$-equivariant line bundles on $\PP$ are of the form $L_{u,v}(i) = L_{u,v}\oo \OO_{\PP}(i)$, and any $T$-equivariant vector bundle splits as a direct sum of $T$-equivariant line bundles (this is proved in \cite{kumar} as a generalization of the classical result by Grothendieck in the non-equivariant setting). Since the bundles $\mc{F}^d_r$ are $T$-equivariant (in fact, they are $\GL_2$-equivariant), they split as direct sums of $L_{u,v}(i)$, and \cite[Theorem~3.2]{KMRR} gives a recursive description of this splitting.

We implement the method \texttt{splittingFdr} to describe the splitting type of the bundles $\mc{F}^d_r$. The method takes inputs $p,d,r$ and returns the splitting type of $\mc{F}^d_r$ over a field of characteristic~$p$. It has an optional input \texttt{Multidegree} having default value \texttt{false}. If \texttt{Multidegree} is set to \texttt{true} then the output describes the $T$-equivariant splitting as a list of triples $\{i,u,v\}$, each corresponding to a summand $L_{u,v}(i)$. If \texttt{Multidegree} is set to \texttt{false} then the output records the non-equivariant splitting (ignoring the $1$-dimensional factors $L_{u,v}$). 

\begin{ex}\label{ex:split-F15-7}
The splitting type of $\mc{F}^{15}_7$ in characteristic $5$ is computed as follows:
\smallskip

\noindent\extt{i1: loadPackage "IncidenceCorrespondenceCohomology"\\
i2: splittingFdr(5,15,7)\\
o2 = \{-10, -8, -10, -8, -9, -9, -9\}\\
o2: List}

\noindent This computation signifies the existence of an isomorphism
\[\mc{F}_7^{15} \simeq \mc{O}_{\PP}(-10)^{\oplus 2}\oplus\mc{O}_{\PP}(-9)^{\oplus 3}\oplus\mc{O}_{\PP}(-8)^{\oplus 2}\]
To compute the $T$-equivariant splitting, we proceed as follows:
\smallskip

\noindent\extt{i3: splittingFdr(5,15,7, Multidegree=>true)\\
o3 = \{\{-10, 15, 10\}, \{-8, 14, 9\}, \{-10, 10, 15\}, \{-8, 9, 14\}, \{-9, 13, 11\},\\
\phantom{o3 =} \{-9, 12, 12\}, \{-9, 11, 13\}\}\\
o3: List}

\noindent This refines the above splitting to a $T$-equivariant isomorphism
\[\mc{F}_7^{15} \simeq L_{15,10}(-10) \oplus L_{10,15}(-10) \oplus L_{13,11}(-9)\oplus L_{12,12}(-9)\oplus L_{11,13}(-9)\oplus L_{14,9}(-8)\oplus L_{9,14}(-8) \]
\end{ex}

\begin{ex} The splitting type of the bundles $\mc{F}^d_r$ can also be computed by direct calculation, considering the minimal generators of the graded module associated to $\mc{F}^d_r$. Below is a time comparison between our implementation based on \cite[Theorem 3.2]{KMRR} and the direct calculation:
\smallskip

      \begin{tabularx}{0.95\textwidth} { 
  | >{\centering\arraybackslash}X 
  | >{\centering\arraybackslash}X | }
 \hline
    \texttt{splittingFdr(5,2249,1112)} &  Direct Calculation \\
 \hline
    0.0007495 seconds  &      188.788 seconds\\
\hline
\end{tabularx}
\smallskip

\noindent For the $T$-equivariant description we have:
\smallskip

      \begin{tabularx}{0.95\textwidth} { 
  | >{\centering\arraybackslash}X 
  | >{\centering\arraybackslash}X | }
 \hline
    \texttt{splittingFdr(5,2249,1112, Multidegree=>true)} &  \vspace{0.1mm}Direct Calculation \\
 \hline
  0.211905 seconds &     180.982 seconds\\
\hline
\end{tabularx}
\end{ex}

\subsection{Principal parts on $\PP$}\label{subsec:pparts}

We write $U=\kk^2 = H^0(\PP,\OO_\PP(1))$ in order to better keep track of equivariance, and note that we have a $T$-equivariant decomposition $U=L_{1,0}\oplus L_{0,1}$. We write $\det U=\bigwedge^2 U \simeq L_{1,1}$, and consider the diagonal embedding $\PP\hookrightarrow\PP\times\PP$, with ideal sheaf $\mc{I} = \det U \oo \OO_{\PP\times\PP}(-1,-1)$. The order $k+1$ thickening of $\PP$ is the subscheme $\PP^{(k)} \subset \PP\times\PP$ defined by $\mc{I}^{k+1}$, and we have a short exact sequence
\begin{equation}\label{eq:ses-pparts} 0 \lra (\det U)^{k+1} \oo \OO_{\PP\times\PP}(-k-1,-k-1) \lra \OO_{\PP\times\PP} \lra \OO_{\PP^{(k)}} \lra 0.
\end{equation}
If we denote by $\pi_1,\pi_2:\PP^{(k)}\lra \PP$ the natural projections, then to any line bundle $\OO_\PP(m)$ we can associate the \textbf{vector bundle of $k$-th order principal parts}
\[ \mc{P}^k(\OO_\PP(m)) = \pi_{1*}\pi_2^*(\OO_\PP(m)) = \pi_{1*}(\OO_{\PP^{(k)}}(0,m)).\]
Twisting \eqref{eq:ses-pparts} by $\OO_{\PP\times\PP}(0,m)$ and pushing forward along $\pi_1$ yields three possibilities. 
\smallskip

\noindent{\bf Case 1.} If $m\geq k+1$ then we get a short exact sequence
\[ 0 \lra (\det U)^{k+1} \oo \Sym^{m-k-1} U \oo \OO_\PP(-k-1) \lra \Sym^m U \oo \OO_\PP \lra \mc{P}^k(\OO_\PP(m)) \lra 0 \]
Dualizing, tensoring by $(\det U)^m$, and using the identification $(\Sym^m U)^{\vee} \oo (\det U)^m = D^m U$, we get an exact sequence
\[ 0 \lra \left(\mc{P}^k(\OO_\PP(m))\right)^{\vee} \oo (\det U)^m \lra D^m U \oo \OO_\PP \lra D^{m-k-1}U \oo \OO_\PP(k+1) \lra 0\]
Letting $d=m$, $r=k+1$, and comparing with \eqref{eq:ses-Fdr} we conclude that
\[ \mc{P}^k(\OO_\PP(m)) = \left(\mc{F}^m_{k+1}\right)^{\vee} \oo (\det U)^m = \left(\mc{F}^m_{k+1}\right)^{\vee} \oo L_{m,m} \]
so the ($T$-equivariant) splitting of $\mc{P}^k(\OO_\PP(m))$ is determined by our earlier calculation.

\smallskip

\noindent{\bf Case 2.} If $-1\leq m\leq k$ then we get a short exact sequence
\[ 0 \lra \Sym^m U \oo \OO_\PP \lra \mc{P}^k(\OO_\PP(m)) \lra (\det U)^k \oo (\Sym^{k-1-m}U)^{\vee} \oo \OO_\PP(-k-1) \lra 0\]
Since $\Ext^1(\OO_\PP(-k-1),\OO_\PP)=\HH^1(\PP,\OO_\PP(k+1))=0$, the above sequence splits and we get
\[ 
\begin{aligned}
    \mc{P}^k(\OO_\PP(m)) &\simeq \left(\bigoplus_{i=0}^m L_{i,m-i}(0) \right)  \oplus \left(\bigoplus_{i=0}^{k-1-m} L_{m+1+i,k-i}(-k-1) \right) \\
    & \simeq \OO_\PP^{\oplus (m+1)} \oplus \OO_\PP(-k-1)^{\oplus (k-m)}
\end{aligned}
\]

\smallskip

\noindent{\bf Case 3.} If $m\leq -2$ then we get a short exact sequence
\[
\begin{aligned}
0 \lra \mc{P}^k(\OO_\PP(m)) &\lra (\det U)^k \oo (\Sym^{k-1-m}U)^{\vee} \oo \OO_\PP(-k-1) \lra \\
&\lra (\det U \oo \Sym^{-m-2}U)^{\vee} \oo \OO_\PP \lra 0
\end{aligned}
\]
Tensoring by $(\det U)^{-1-m} \oo \OO_\PP(k+1)$ and passing to divided powers as in {\bf Case 1} we get a short exact sequence
\[0 \lra \mc{P}^k(\OO_\PP(m)) \oo (\det U)^{-1-m} \oo \OO_\PP(k+1) \lra D^{k-1-m}U \oo \OO_\PP \lra D^{-2-m}U \oo \OO_\PP(k+1) \lra 0 \]
Letting $d=k-1-m$, $r=k+1$, and comparing with \eqref{eq:ses-Fdr} we conclude that
\[ \mc{P}^k(\OO_\PP(m)) = \mc{F}^{k-1-m}_{k+1} \oo (\det U)^{1+m} \oo \OO_\PP(-k-1)\simeq \mc{F}^{k-1-m}_{k+1} \oo L_{1+m,1+m}(-k-1) \]
and the ($T$-equivariant) splitting of $\mc{P}^k(\OO_\PP(m))$ is again determined by that of the bundles~$\mc{F}^d_r$.

Having implemented a function that computes the (T-equivariant) splitting type of $\mc{F}^d_r$, we also implement a method \texttt{splittingPrincipalParts} based on the relationship between $\mc{F}_r^d$ and $\mc{P}^k(\OO_\PP(m))$ discussed above. The method has as inputs $p, m, k$ and returns the splitting type of $\mc{P}^k(\OO_\PP(m))$ over a field of characteristic $p$. Like the function \texttt{splittingFdr}, \texttt{splittingPrincipalParts} also has the optional Boolean input of \texttt{Multidegree} which is defaulted to \texttt{false}. The output of \texttt{splittingPrincipalParts} is also formatted in the same manner as \texttt{splittingFdr}. If \texttt{Multidegree} is set to \texttt{true}, then the output is a list of triples \texttt{\{i,u,v\}} corresponding to summands $L_{u,v}(i)$ that occur in the decomposition. In the default case when \texttt{Multidegree} is false, the output is only a list of integers corresponding to the summands $\OO_{\PP}(i)$ that appear.

\begin{ex} We compute the splitting type of $\mc{P}^6(\OO_{\PP}(15))$ in characteristic $5$ (see also Example~\ref{ex:split-F15-7}):
\smallskip

\noindent\extt{i4: splittingPrincipalParts(5,15,6)\\
o4 = \{10, 8, 10, 8, 9, 9, 9\}\\
o4: List}

\noindent\extt{i5: splittingPrincipalParts(5,15,6, Multidegree=>true)\\
o5 =  \{\{10, 0, 5\}, \{8, 1, 6\}, \{10, 5, 0\}, \{8, 6, 1\}, \{9, 2, 4\}, \{9, 3, 3\}, \{9, 4, 2\}\}
o5: List}
\end{ex}

\section{Graded Han-Monsky representation ring}\label{sec:han-monsky}

In this section we discuss the multiplication in the \textbf{graded Han--Monsky representation ring} (see \cite{han-monsky} and \cite[Section~5]{KMRR}), which is the Grothendieck ring of the category of finite length graded $\kk[T]$-modules, with tensor product $M\oo_{\kk} N$ defined by
\[ T \cdot (m\oo n) = Tm\oo n + m\oo Tn\quad\text{ for }m\in M,\ n\in N.\]
The indecomposable graded $\kk[T]$-modules are of the form $\kk[T]/(T^d)(-j)$ where $d\geq 1$ denotes the length of the module and $j\in\ZZ$ denotes the degree of the cyclic generator. We write $\delta_d(-j)$ for the corresponding isomorphism class. For $0\leq j<a\leq b$, there exist unique non-negative integers $c_j = c_j(a,b)$ such that   
\begin{equation}\label{eq:dela-delb-general} \delta_a\cdot \delta_b = \sum_{j=0}^{a-1}\delta_{c_j}(-j),
\end{equation}
and understanding the multiplication in the Han--Monsky ring amounts to identifying $c_j(a,b)$. In characteristic zero one has $c_j(a,b)=a+b-2j-1$, but such an explicit formula in positive characteristic remains unknown. We discuss the recursive calculation of $c_j(a,b)$ in Section~\ref{subsec:rec-cjab}.

Of particular interest is the calculation of the $n$-fold product $\delta_{a_1}\cdots\delta_{a_n}$, which represents the Artinian monomial complete intersection
$$A = \kk[T_1,\cdots,T_n]/\langle T_1^{a_1},\cdots,T_n^{a_n}\rangle$$
viewed as a $\kk[T]$-module by letting $T=T_1+\cdots+T_n$. An explicit formula for $\delta_{a_1}\cdots\delta_{a_n}$ would describe the Jordan type of $A$ with respect to the linear form $T_1+\cdots+T_n$ (for more about Jordan type, see \cite{AIM}), which in turn would determine the (weak and strong) Lefschetz properties for~$A$. We discuss this in more detail in Section~\ref{subsec:lefschetz}.

\subsection{Recursive description of the multiplication}\label{subsec:rec-cjab}

We implement the method \texttt{hanMonsky} to compute the product $\delta_{a_1}\cdots\delta_{a_n}$ over a field $\kk$. The method takes as inputs the characteristic~$p$ of~$\kk$, and a list $L=\{a_1,\cdots,a_n\}$. It outputs a HashTable with entries of the form 
      $$ \{c \Rightarrow f(q)\}$$ 
where $c\geq 1$ and $f(q)\in\ZZ[q,q^{-1}]$ is a Laurent polynomial with non-negative coefficients. The polynomial $f(q)$ is uniquely determined by the condition that
$$ f_j\cdot q^j \text{ is a term in } f(q) \iff \delta_c(-j) \text{ is a summand of  } \delta_{a_1}\cdot\delta_{a_2}\cdots\delta_{a_n} \text{ with multiplicity }f_j.$$
Our method for computing the product is recursive, with base case the $2$-fold products $\delta_a\delta_b$. The method \texttt{hanMonsky} also takes an optional input \texttt{UseConjecture} with default value \texttt{true}, in which case the products $\delta_a\cdot\delta_b$ are computed using the following conjectural recursive description.

\begin{con}\label{Conjecture}
    Given $0\leq j<a\leq b$, we can give a recursive description of the integers $c_j$ in \eqref{eq:dela-delb-general} as follows. 
   Let $e\geq 0$ such that $q'=p^{e-1}<a\leq q=p^e$, and let $r$ such that
    \[ rq \leq a+b-1-j < (r+1)q.\]
    \begin{itemize}
        \item If $b-j\leq rq$ then
   $ c_j(a,b) = rq.$
    \end{itemize}
 Otherwise, (when $b-j>rq$) we define $m$ so that
   $mq'\leq a<(m+1)q'$.    
    Let $a'=a-mq'$ and consider $i$ such that
    $ iq'\leq j\leq (i+1)q'-1.$
    We have
    \begin{itemize}
        \item If $j\leq iq'+a'-1$ then 
        $c_j(a,b) = c_{j-iq'}(a',b+(m-2i)q').$ 
        \item If $j\geq iq'+a'$ then 
       $c_j(a,b) = c_{j-iq'-a'}(q'-a',b+(m-1-2i)q').$
    \end{itemize}
\end{con}

When the option \texttt{UseConjecture} is set to \texttt{false}, the products $\delta_a\cdot\delta_b$ are computed by a direct algebraic method, which is slower than the alternative above.

\begin{ex} We compute the product  $\delta_4\delta_6$ in characteristic $3$, first using Conjecture~\ref{Conjecture}:
\smallskip

\noindent\extt{i1: loadPackage "IncidenceCorrespondenceCohomology"\\
i2: hanMonsky(3,\{4,6\})\\
o2 =  HashTable$\mathtt{\left\{3\ \Rightarrow \ q^{3},\,6\ \Rightarrow \ q^{2}+q,\,9\ \Rightarrow \ 1\right\}}$\\
o2: HashTable\\}
so we have \[\delta_4\delta_6=\delta_9+\delta_6(-1)+\delta_6(-2)+\delta_3(-3).\]
Next, we recompute this without using Conjecture~\ref{Conjecture}. \smallskip

\noindent\extt{i3: hanMonsky(3,\{4,6\}, UseConjecture=>false)\\
o3 =  HashTable$\mathtt{\left\{3\ \Rightarrow \ q^{3},\,6\ \Rightarrow \ q^{2}+q,\,9\ \Rightarrow \ 1\right\}}$\\
o3: HashTable}
\end{ex}

\begin{ex}We compute the product $\delta_3\delta_4\delta_6$ in different characteristics. 
\smallskip

\noindent\extt{i4: hanMonsky(2,\{3,4,6\})\\
o4 =  HashTable$\mathtt{\left\{4\ \Rightarrow \ q^5  + 2q^4  + 2q^3  + q^2,\,8\ \Rightarrow \ q^3  + 2q^2  + 2q + 1\right\}}$\\
o4: HashTable\\
i5: hanMonsky(3,\{3,4,6\})\\
o5 = HashTable$\mathtt{\left\{3\ \Rightarrow \ q^{5}+q^{4}+q^{3},\,6\ \Rightarrow \ q^{4}+2\,q^{3}+2\,q^{2}+q,\,9\ \Rightarrow \ q^{2}+q+1\right\}}$\\
o5: HashTable\\
i6: hanMonsky(5,\{3,4,6\})\\
o6 =  HashTable$\mathtt{\left\{5\ \Rightarrow \ q^{5}+2\,q^{4}+3\,q^{3}+2\,q^{2}+q,\,7\ \Rightarrow \ q^{2},\,10\ \Rightarrow \ q+1\right\}}$\\
o6: HashTable\\
i7: hanMonsky(7,\{3,4,6\})\\
o7 = HashTable$\mathtt{\left\{1\ \Rightarrow \ q^{5},\,3\ \Rightarrow \ q^{4},\,5\ \Rightarrow \ q^{3},\,7\ \Rightarrow \ q^{4}+2\,q^{3}+3\,q^{2}+2\,q+1\right\}}$\\
o7: HashTable}

\noindent For $p>10$ the product is the same as in characteristic zero, given by:
\smallskip

\noindent\extt{o8: HashTable$\mathtt{\left\{1\ \Rightarrow \ q^{5},\,3\ \Rightarrow \ 2\,q^{4},\,5\ \Rightarrow \ 3\,q^{3},\,7\ \Rightarrow \ 3\,q^{2},\,9\ \Rightarrow \ 2\,q,\,11\ \Rightarrow \ 1\right\}}$
}
\end{ex} 

We compare below the running times for some examples with the option \texttt{UseConjecture} set as \texttt{true} (default setting) and set as \texttt{false}.
\begin{ex} The table below compares the times to compute the Han--Monsky multiplication of $\delta_3\delta_8\delta_{14}\delta_{31}$ in characteristic $3$ using the conjecture, \texttt{hanMonsky(3,\{3,8,14,31\})}, and without it, \texttt{hanMonsky(3,\{3,8,14,31\}, UseConjecture=>false)}

\smallskip

      \begin{tabularx}{0.95\textwidth} {
  | >{\centering\arraybackslash}X 
  | >{\centering\arraybackslash}X | }
 \hline
   \texttt{hanMonsky(3,\{3,8,14,31\})} &  \texttt{UseConjecture=>false} \\
 \hline
 0.0015699  seconds&  0.462372 seconds\\
\hline
\end{tabularx}

\end{ex}
\begin{ex} Time comparison for the computation of $\delta_{38}\delta_{14}\delta_{51}$ in $\ch\kk=7$ 
\smallskip

      \begin{tabularx}{0.95\textwidth} { 
  | >{\centering\arraybackslash}X 
  | >{\centering\arraybackslash}X | }
 \hline
 \texttt{hanMonsky(7,\{38,14,51\})}  &  \texttt{UseConjecture=>false} \\
 \hline
0.0020716 seconds &  0.311123 seconds\\
\hline
\end{tabularx}

\end{ex}
\begin{ex} Time comparison for the computation of $\delta_3\delta_7\delta_8\delta_{14}\delta_{21}$ in $\ch\kk=5$ 
\smallskip

      \begin{tabularx}{0.95\textwidth} { 
  | >{\centering\arraybackslash}X 
  | >{\centering\arraybackslash}X | }
 \hline
\texttt{hanMonsky(5,\{3,7,8,14,21\})} &  \texttt{UseConjecture=>false} \\
 \hline
 0.0027694 seconds &  0.363145 seconds\\
\hline
\end{tabularx}

\end{ex}

\medskip

\subsection{Lefschetz Properties}\label{subsec:lefschetz} 

Recall that an Artinian algebra $A$ has the Weak Lefschetz Property (WLP) if there exists a linear form $\ell\in A_1$ such that the multiplication maps $A_i \overset{\times\ell}{\to} A_{i+1}$ have maximal rank for all $i$, and it has the Strong Lefschetz Property (SLP) if $A_i \overset{\times\ell^d}{\to} A_{i+d}$ has maximal rank for all $i,d$. Here we are interested in monomial complete intersections $A = \kk[T_1,\cdots,T_n]/\langle T_1^{a_1},\cdots,T_n^{a_n}\rangle$, in which case the Lefschetz properties can be tested on the linear form $\ell=T_1+\dots T_n$ \cite{MMN}.  It is well-known that in $\ch{\kk}=0$, $A$ has both the SLP and WLP  \cite{stanley,watanabe}  but the problem is more subtle in positive characteristic. The monomial complete intersections having the SLP are classified by \cite[Theorem 3.8]{lund-nick} and \cite[\S 3]{nick} based on different techniques, while a similar classification for the WLP is still an open question. 

Knowing the multiplication in the Han--Monsky ring allows for easy criteria to test WLP and SLP: if we let $s=a_1+\cdots+a_n-n$ denote the socle degree of $A$ then
\begin{enumerate}
    \item WLP holds if and only if every summand $\delta_c(-j)$ of $\delta_{a_1}\cdots\delta_{a_n}$ satisfies $ j + (c-1) \geq \frac{s}{2}$; 
    \item SLP holds if and only if every summand $\delta_c(-j)$ of $\delta_{a_1}\cdots\delta_{a_n}$ satisfies $2j + (c-1) = s$.
\end{enumerate}

\begin{ex}\label{ex:346-p=3} Using the above criteria, we can check that $A=\kk[T_1,T_2,T_3]/\langle T_1^3,T_2^4,T_3^6\rangle$ satisfies WLP but fails SLP when $\ch\kk=3$:
\smallskip

\noindent\extt{i9: HM = hanMonsky(3,\{3,4,6\})\\
o9 = HashTable$\mathtt{\left\{3\ \Rightarrow \ q^{5}+q^{4}+q^{3},\,6\ \Rightarrow \ q^{4}+2\,q^{3}+2\,q^{2}+q,\,9\ \Rightarrow \
      q^{2}+q+1\right\}}$\\
o9: HashTable}
\end{ex}
 
We implement the methods \texttt{hasWLP} and \texttt{hasSLP} to determine if a monomial complete intersection satisfies WLP and SLP respectively. The methods take as inputs the characteristic $p$ of the field $\kk$ and a list $L=\{a_1,\dots,a_n\}$ of exponents describing the monomial complete intersection. The output is the Boolean value \texttt{true} when the relevant Lefschetz property holds, and it is \texttt{false} otherwise. Both \texttt{hasWLP} and \texttt{hasSLP} take optional input \texttt{UseConjecture}, which refers as before to the use of Conjecture~\ref{Conjecture} for computing products in the Han--Monsky ring.

\begin{ex} We consider the monomial complete intersections  $$A=\kk[T_1,T_2,T_3,T_4]/\langle T_1^{3},T_2^4,T_3^6,T_4^{8}\rangle \quad \text{and} \quad A'=\kk[T_1,T_2,T_3,T_4]/\langle T_1^{3},T_2^4,T_3^6,T_4^{13}\rangle$$ in characteristic $7$:
\smallskip
 
\noindent\extt{i10: hasWLP(7,\{3,4,6,8\})\\
o10 = false\\
i11: hasWLP(7,\{3,4,6,13\})\\
o11 = true\\
i12: hasSLP(7,\{3,4,6,13\})\\
o12 = false}
\smallskip

\noindent This shows that $A$ fails WLP (and therefore also SLP), while $A'$ satisfies WLP but not SLP.
\end{ex}

The default implementation for \texttt{hasSLP} has \texttt{UseConjecture} set to \texttt{false}, and is based on the criteria from \cite{nick,lund-nick}. When \texttt{UseConjecture} is set to \texttt{true}, we combine Conjecture~\ref{Conjecture} with the following reformulation of criterion (2) above: SLP holds if and only if all the partial products $\delta_{a_1}$, $\delta_{a_1}\delta_{a_2}$, $\cdots$, $\delta_{a_1}\cdots\delta_{a_n}$ are the same as in characteristic zero.

For the method \texttt{hasWLP}, we return \texttt{true} if $\op{char}(\kk)=0$ or $n\leq 2$, and we apply the criterion from \cite[Theorem~8.1]{KMRR} if $\op{char}(\kk)=2$, ignoring the option \texttt{UseConjecture}. For all other cases, if \texttt{UseConjecture} is set to \texttt{true} (the default value) then we use criterion (1) above together with the conjectural Han--Monsky multiplication to test WLP. If \texttt{UseConjecture} is set to \texttt{false}, then we apply \cite[Proposition 3.5(2)]{HMMNWW}, which implies that WLP is equivalent to 
\[\dim_{\kk}(A/TA) = \dim_{\kk}A_{\lfloor s/2\rfloor}. \]
The right side of the above equality measures the maximal value of the Hilbert function of $A$, also called the \textbf{Sperner number of $A$}, and it is independent of $\op{char}(\kk)$. The left side computes the minimal number of generators of $A$ as a $\kk[T]$-module, or equivalently, the number of terms in the expansion of the product $\delta_{a_1}\cdots\delta_{a_n}$, counted with multiplicity. In the language of \cite{AIM}, this is also the number of parts for the Jordan type of $A$ (which records the summands $\delta_c$ of $\delta_{a_1}\cdots\delta_{a_n}$ with multiplicity, ignoring their degree shift).

\begin{ex} With the notation from Example~\ref{ex:346-p=3}, we can compute the Jordan type of $A$ by evaluating each of the Laurent polynomials in the HashTable \texttt{HM} at $q=1$:
\smallskip

\noindent\extt{i13: jordanType = flatten apply(keys HM,c -> splice\{HM\#c[1] : c\}) \\
o13 = \{9,9,9,6,6,6,6,6,6,3,3,3\}\\
o13: List\\
i14: \#jordanType\\
o14 = 12\\
}
\smallskip
We can perform similar calculations in characteristic zero:
\smallskip

\noindent\extt{i15: HM0 = hanMonsky(0,\{3,4,6\})\\
o15 = $\mathtt{HashTable}\left\{1\ \Rightarrow \ q^{5},\:3\ \Rightarrow \ 2\,q^{4},\:5\ \Rightarrow \ 3\,q^{3},\:7\ \Rightarrow \ 3\,q^{2},\:9\ \Rightarrow \ 2\,q,\:11\ \Rightarrow \ 1\right\}$\\
o15: HashTable\\
i16: jordanType0 = flatten apply(keys HM0,c -> splice\{HM0\#c[1] : c\}) \\
o16 = \{1, 9, 9, 3, 3, 11, 5, 5, 5, 7, 7, 7\}\\
o16: List\\
i17: \#jordanType0\\
o17 = 12\\}
Notice that although the decomposition of $\delta_3\delta_4\delta_6$ is very different in characteristics $0$ and $3$, the number of components (or Jordan blocks) is the same, which is what characterizes WLP. On the other hand, the fact that \texttt{HM} and \texttt{HM0} are different explains the failure of SLP in characteristic $3$.
\end{ex}

As an application of the method \texttt{hasWLP}, we implement the method \texttt{monomialCIsWithoutWLP} to generate all the monomial complete intersections with fixed embedding dimension and socle degree. The inputs are the characteristic $p$ of $\kk$, the number of variables $n$, and the socle degree~$s$. The output is a list of $n$-tuples $\{a_1,\dots,a_n\}$ with $2\leq a_1\leq \dots \leq a_n$ and $a_1+\dots +a_n-n=s$ such that the corresponding monomial complete intersections fail WLP. As usual, we allow the optional input \texttt{UseConjecture}, which refers to the use of Conjecture~\ref{Conjecture}.

\begin{ex}\label{ex:moncis-fail-wlp}
We compute the monomial complete intersections of embedding dimension $4$ and socle degree $10$ that fail WLP in characteristic $5$.
\smallskip

\noindent\extt{i18: monomialCIsWithoutWLP(5,4,10)\\
o18 = \{\{2, 2, 5, 5\}, \{2, 3, 4, 5\}, \{2, 4, 4, 4\}, \{3, 3, 3, 5\}, \{3, 3, 4, 4\}\}}
\smallskip

\noindent This corresponds to the following algebras that fail WLP in characteristic $5$:
$$\frac{\kk[T_1,\dots,T_4]}{\langle T_1^2,T_2^2,T_3^5,T_4^5\rangle} \quad  \frac{\kk[T_1,\dots,T_4]}{\langle T_1^2,T_2^3,T_3^4,T_4^5\rangle} \quad \frac{\kk[T_1,\dots,T_4]}{\langle T_1^2,T_2^4,T_3^4,T_4^4\rangle} \quad \frac{\kk[T_1,\dots,T_4]}{\langle T_1^3,T_2^3,T_3^3,T_4^5\rangle}  \quad \frac{\kk[T_1,\dots,T_4]}{\langle T_1^3,T_2^3,T_3^4,T_4^4\rangle} $$
\end{ex}
We compare below the running times for the method \texttt{monomialCIsWithoutWLP} with the option \texttt{UseConjecture} set as \texttt{true} (default setting) and set as \texttt{false}.
\begin{ex} Time comparison for \texttt{monomialCIsWithoutWLP} as in Example~\ref{ex:moncis-fail-wlp}:
\smallskip

      \begin{tabularx}{0.95\textwidth} {
  | >{\centering\arraybackslash}X 
  | >{\centering\arraybackslash}X | }
 \hline
   \texttt{monomialCIsWithoutWLP(5,4,10)} &  \texttt{UseConjecture=>false} \\
 \hline
 0.0020031 seconds &  0.024318 seconds\\
\hline
\end{tabularx}

\end{ex}
\begin{ex} Time comparison for \texttt{monomialCIsWithoutWLP} for $n=6$ and socle degree  $s=30$ in characteristic $7$:
\smallskip

 \begin{tabularx}{0.95\textwidth} {
  | >{\centering\arraybackslash}X 
  | >{\centering\arraybackslash}X | }
 \hline
   \texttt{monomialCIsWithoutWLP(7,6,30)} &  \texttt{UseConjecture=>false} \\
 \hline
  1.25439 seconds &  5.77864 seconds \\
\hline
\end{tabularx}
\end{ex}

We also implement an instance of \texttt{hasWLP} that can be used to check WLP for a graded Artinian algebra $R/I$, where $R=\kk[x_1,\dots, x_n]$ is a polynomial ring over a sufficiently large field. This method takes optional inputs \texttt{GorensteinAlg} and \texttt{MonomialAlg} that allow faster computations for Gorestein and monomial ideals respectively. 

 The instance \texttt{hasWLP(R,I)} requires $R$ to be a standard graded polynomial ring over a sufficiently large field $\kk$, e.g., $\mathtt{R=QQ[x_1,…,x_n]}$. When $R$ is a polynomial ring over a finite field, then when \texttt{hasWLP(R, I)} outputs \texttt{false}, it confirms that $R/I$ fails the WLP. 
   However, when \texttt{hasWLP(R, I)} outputs \texttt{true}, we can only conclude that $R/I$ has the WLP over a field extension. In the case that~$I$ is a monomial ideal, this is no longer a concern and the test works over any field.

\begin{ex}The following is an example of a Gorenstein Algebra that fails the WLP: \smallskip

\noindent\extt{i19: $\mathtt{R=QQ[x,y,z,w,t]}$;\\
i20:  $\mathtt{F=x^4yzt +x^2y^2t^2w}$;\\
i21: I=inverseSystem(F)\\
o21 = ideal $\mathtt{(w^2 , zw, z^2 , t^3 , zt^2 , y^2 z, x^2 z - 3ywt, y^3 , x^3 w, x^3 t^2 , x^3 y^2 , x^5 )}$\\
o21: Ideal of R\\
i22: hasWLP(R,I)\\
o22 = false}\\
Since in this example we know that $R/I$ is Gorenstein, we can use the option \texttt{GorensteinAlg}:\\
\extt{i23: hasWLP(R,I, GorensteinAlg => true)\\
o23 = false}\smallskip

\noindent We note that a simple sufficient condition to guarantee WLP when $F$ is a binomial is given in \cite[Theorem 3.4]{BinMaculayDual}, and a construction of examples of Gorenstein algebras that fail WLP in codimension $\geq 4$ appears in \cite[Example~3.7]{BinMaculayDual}.
\end{ex}

\begin{ex} The almost complete intersection $A=\mathbb{Q}[x,y,z]/\langle x^9,y^9,z^9, x^3y^3z^3\rangle$ fails~WLP: 
\noindent\extt{i24: $\mathtt{R=QQ[x,y,z]}$;\\
i25: I=ideal($\mathtt{x^9,y^9,z^9,x^3y^3z^3}$)\\
i26:  hasWLP(R,I,GorensteinAlg => false)\\
o26 = false}\\
Since in this example $I$ is a monomial, we can also use the option \texttt{MonomialAlg}: \\
\texttt{i27: hasWLP(R,I, MonomialAlg =>true)\\
o27 = false}\\
\end{ex}

\section*{Acknowledgements}
 The authors would like to thank Mike Stillman,  Keller VandeBogert, and Matthew Weaver for helpful discussions regarding various aspects of this project. Marangone gratefully acknowledges that this research was supported in part by the Pacific Institute for the Mathematical Sciences. 
 Raicu and Reed acknowledge the support of the National Science Foundation Grant DMS-2302341. Part of the material in this paper is based upon work supported by the National Science Foundation under Grant No. DMS-1928930 and by the Alfred P. Sloan Foundation under grant G-2021-16778, while Raicu and Reed were in residence at the Simons Laufer Mathematical Sciences Institute (formerly MSRI) in Berkeley, California, during the Spring 2024 semester. Part of the work on this project was done while Marangone was in residence at the Fields Institute, Toronto, during the \emph{Thematic Program in Commutative Algebra and Applications}, during the Winter 2025 semester.

\bibliographystyle{amsalpha}
\bibliography{IncidenceCorrespondenceCohomologyPackage}

\end{document}